\documentclass{article}
\usepackage{amsfonts}
\usepackage{amssymb}
\usepackage{amsmath}
\usepackage{srcltx}
\usepackage{inputenc}
\inputencoding{ansinew}

\newcommand{\bnom}{\begin{nom}}
\newcommand{\enom}{\end{nom}}
\newcommand{\theo}{{\rm\bf Theorem}}

\newcommand{\defi}{{\rm\bf Definition}}
\newcommand{\ob}{{\rm\bf Remark}}

\newcommand{\corol}{{\rm\bf Corollary}}

\newcommand{\demo}{{\bf Proof}}

\newcommand{\ex}{{\bf Example}}
\newcommand{\qed}{\hfill\mbox{$_{_{\blacksquare}}$}}

\newtheorem{nom}{{\!}}[section]
\addtocounter{nom}{0}

\begin{document}
\title{COMMON FIXED POINT THEOREMS FOR OCCASIONALLY WEAKLY COMPATIBLE MAPS}
\author{ H. Bouhadjera {\footnotesize {\dag} } \& C. Godet-Thobie
{\footnotesize {\ddag}} \\
{\normalsize 1-Universit\'{e} europ\'{e}enne de Bretagne, France.}\\
 {\normalsize
2-Universit\'{e} de Bretagne Occidentale,}\\
{\normalsize Laboratoire de Math\'{e}matiques de Brest; Unit\'{e}
CNRS: UMR
6205 }\\
{\normalsize 6, avenue Victor Le Gorgeu, CS 93837, F-29238 BREST Cedex 3
FRANCE }}
\date{ {\footnotesize E-Mail: \dag hakima.bouhadjera@univ-brest.fr; \ddag
christiane.godet-thobie@univ-brest.fr } }
\maketitle

\begin{abstract}
In this paper, we establish some common fixed point theorems for two
pairs of occasionally weakly compatible single and set-valued maps
satisfying a strict contractive condition in a metric space. Our
results unify and extend many results existing in
the literature including those of Aliouche \cite{Ali}, Bouhadjera
\cite{Bou} and Popa \cite{Pop1}-\cite{Pop6}. Also we establish another
common fixed point theorem for four occasionally weakly compatible single and set-valued maps of Gregu\v{s} type which improves the results of
Djoudi and Nisse \cite{DjN}, Pathak et al. \cite{PCKM} and others
and we end our work by giving another theorem which generalizes
the results given by Elamrani and Mehdaoui \cite{ElM}, Mbarki
\cite{Mba} and references therein.


\end{abstract}

\section{Introduction and preliminaries}

Throughout this paper, $(\mathcal{X},d)$ denotes a metric space and $\mathcal{P}_{fb}(\mathcal{X})$ the class of all nonempty bounded closed subsets of $\mathcal{X}$. We recall these usual notations:  for $x\in \mathcal{X}$ and $A\subseteq \mathcal{X}$,
\begin{eqnarray*}
d(x,A)=\inf \{d(x,y):y\in A\}.
\end{eqnarray*}
Let $H$ be the associated Hausdorff metric on $\mathcal{P}_{fb}(\mathcal{X})$: for every $A$ and every $B$ in $\mathcal{P}_{fb}(\mathcal{X})$,
\begin{eqnarray*}
H(A,B)=\max \{\sup_{x\in A}d(x,B),\sup_{y\in B}d(A,y)\}
\end{eqnarray*}
and
\begin{eqnarray*}
\delta (A,B)=\sup \{d(a,b):a\in A,b\in B\}.
\end{eqnarray*}
For simplicity, we write $\delta(a,B)$ in place of $\delta(\{a\},B)$; as well as $\delta(A,b)$ in place of $\delta(A,\{b\})$.

In the following, we use small letters: $f$, $g$, $\ldots$ to denote maps
from $\mathcal{X}$ to $\mathcal{X}$ and capital letters: $F$, $G$, $\ldots$ for set-valued maps; that is, maps from $\mathcal{X}$ to
$\mathcal{P}_{fb}(\mathcal{X})$ and we write $fx$ for $f(x)$ and $Fx$ for $F(x)$.

The concepts of weak commutativity, compatibility,
noncompatibility and weak compatibility were frequently used to
prove existence theorems in fixed and common fixed points for
single and set-valued maps satisfying certain conditions in
different spaces. The study of common fixed points on occasionally
weakly compatible maps is new and also interesting.  This notion
which is defined by Al-Thagafi and Shahzad \cite{AlS} and which
is published in 2008, has been used by Jungck and Rhoades
\cite{JuR2} in 2006 and by Abbas and Rhoades \cite{AbR} in 2007.

We begin by a short historic of these different
notions. Generalizing the concept of commuting maps, Sessa \cite{Ses} introduced the concept of weakly commuting maps.
$f$ and $g$ are weakly
commuting if
\begin{eqnarray*}
d(fgx,gfx)\leq d(gx,fx)
\end{eqnarray*}
for all $x\in \mathcal{X}$, where $f$ and $g$ are two self-maps of $(\mathcal{X},d)$.

In 1986, Jungck \cite{Jun1} made more generalized commuting
and weakly commuting
maps called compatible maps. $f$ and $g$ are said to be compatible if
\begin{eqnarray*}
(1)\textit{ \ }\lim_{n\rightarrow \infty }d(fgx_{n},gfx_{n})=0
\end{eqnarray*}
whenever $(x_{n})_{n\in \mathbb{N}} $ is a sequence in
$\mathcal{X}$ such that $\underset{n\rightarrow \infty }{\lim
}fx_{n}=\underset{n\rightarrow \infty }{ \lim }gx_{n}=t$ for some
$t\in \mathcal{X}$. This concept has been useful as a tool for
obtaining more comprehensive fixed point theorems. Clearly,
commuting maps are weakly commuting and weakly commuting maps are
compatible, but neither implication is reversible (see \cite{Jun1}).

Further, the same author with Murthy and Cho \cite{JMC} gave
another generalization of weakly commuting maps by introducing the
concept of compatible maps of type $(A)$. $f$ and $g$ are said to
be compatible of type $(A)$ if
in place of (1) we have the two equalities
\begin{eqnarray*}
\lim_{n\rightarrow \infty }d(fgx_{n},g^{2}x_{n})=0\text{ and }
\lim_{n\rightarrow \infty }d(gfx_{n},f^{2}x_{n})=0.
\end{eqnarray*}
Obviously, weakly commuting maps are compatible of type $(A)$. From
\cite{JMC} it follows that the implication is not reversible.

In their paper \cite{PaK}, Pathak and Khan extended type $(A)$ maps by introducing the concept of compatible maps of
type $(B)$ and compared these maps with compatible and compatible maps of type $(A)$ in normed spaces. To be compatible of type $(B)$, $f$ and $g$ above have to satisfy, in lieu of condition $(1)$, the inequalities
\begin{eqnarray*}
\underset{n\rightarrow \infty }{\lim }d(fgx_{n},g^{2}x_{n})
&\leq &\dfrac{1}{2}\left[\underset{n\rightarrow \infty }{\lim }d(fgx_{n},ft)+\underset{n\rightarrow \infty }{\lim }d(ft,f^{2}x_{n}) \right] \\
&&\text{and} \\
\underset{n\rightarrow \infty }{\lim }d(gfx_{n},f^{2}x_{n})
&\leq &\dfrac{1}{2}\left[\underset{n\rightarrow \infty }{\lim }d(gfx_{n},gt)+\underset{n\rightarrow \infty }{\lim }d(gt,g^{2}x_{n}) \right].
\end{eqnarray*}
It is clear that compatible maps of type $(A)$ are
compatible of type $(B)$. The converse is not true (\cite{PaK}).

In 1998, Pathak et al. \cite{PCKM} introduced an extension of compatibility of type $(A)$ by giving the notion of compatible maps of type $(C)$. $f$ and $g$ are compatible of type $(C)$ if they satisfy the two inequalities
\begin{eqnarray*}
\lim_{n\rightarrow \infty }d(fgx_{n},g^{2}x_{n}) &\leq &\frac{1}{3}\left[
\lim_{n\rightarrow \infty }d(fgx_{n},ft)+\lim_{n\rightarrow \infty
}d(ft,f^{2}x_{n})\right. \\
&&\left. +\lim_{n\rightarrow \infty }d(ft,g^{2}x_{n})\right] \\
\lim_{n\rightarrow \infty }d(gfx_{n},f^{2}x_{n}) &\leq &\frac{1}{3}\left[
\lim_{n\rightarrow \infty }d(gfx_{n},gt)+\lim_{n\rightarrow \infty
}d(gt,g^{2}x_{n})\right. \\
&&\left. +\lim_{n\rightarrow \infty }d(gt,f^{2}x_{n})\right].
\end{eqnarray*}
The same authors gave some examples to show that compatible maps of type $(C)$ need not be neither compatible nor compatible of type $(A)$ (resp., type $(B)$).

In \cite{PCKL} the concept of compatible maps of type $(P)$ was
introduced and compared with compatible and compatible maps of
type $(A)$. $f$ and $g$ are
compatible of type $(P)$ if in lieu of $(1)$ we have
\begin{eqnarray*}
\lim_{n\rightarrow \infty }d(f^{2}x_{n},g^{2}x_{n})=0.
\end{eqnarray*}
Note that compatibility, compatibility of type $(A)$ (resp. $(B)$, $(C)$ and $(P)$) are equivalent if $f$ and $g$ are continuous.

Afterwards, Jungck \cite{Jun2} generalized the compatibility, the
compatibility of type $(A)$, $(B)$, $(C)$ and $(P)$ by introducing
the concept of
weak compatibility. He defines $f$ and $g$ to be weakly compatible if $ft=gt$, $t\in \mathcal{X}$ implies $fgt=gft$.

It is known that all of the above compatibility notions imply
weakly compatible notion, however, there exist weakly compatible
maps which are neither compatible nor compatible of type $(A)$,  $(B)$, $(C)$ and $(P)$ (see \cite{Ali}).

Recently in a paper submitted before 2006 but published only in
2008, Al-Thagafi and Shahzad \cite{AlS}  weakened the concept of
weakly compatible maps by giving the new concept of occasionally
weakly compatible maps. Two self-maps $f$ and $g$ of $\mathcal{X}$
are called occasionally weakly compatible maps (shortly owc) if
there is a point $x$ in $\mathcal{X}$ such that $fx=gx$ at which
$f$ and $g$ commute. This notion is used in 2006 by Jungck and
Rhoades \cite{JuR2} to prove some common fixed point theorems in
symmetric spaces.

In their paper \cite{KaS}, Kaneko and Sessa extended the
compatibility to the
setting of single and set-valued maps as follows: $f:\mathcal{X}
\rightarrow \mathcal{X}$ and $F:\mathcal{X}\rightarrow \mathcal{P}_{fb}(\mathcal{X})$ are
said to be compatible if $fFx\in \mathcal{P}_{fb}(\mathcal{X})$ for all $x\in \mathcal{X}$ and
\begin{eqnarray*}
\lim_{n\rightarrow \infty }H(Ffx_{n},fFx_{n})=0
\end{eqnarray*}
whenever $(x_{n})_{n \in \mathbb{N}}$ is a sequence in $\mathcal{X}$ such that $fx_{n}\rightarrow t$, $Fx_{n}\rightarrow A\in \mathcal{P}_{fb}(\mathcal{X})$ and $ t\in A$.

After, in \cite{JuR1} Jungck and Rhoades extend the
concept of compatible single and set-valued maps by giving the
concept of weak compatibility. Maps $f$ and $F$ are weakly
compatible if they commute at their coincidence points; i.e., if
$fFx=Ffx$ whenever $fx\in Fx$.

More recently, Abbas and Rhoades \cite{AbR} extended the
definition of owc maps to the setting of set-valued maps and they
proved some common fixed point theorems satisfying generalized
contractive condition of integral type. $f$ and $F$ are said to be
owc if and only if there exists some point $x$ in $\mathcal{X}$
such that $fx\in Fx$ and $fFx\subseteq Ffx$. Clearly, weakly
compatible maps are occasionally weakly compatible. However, the
converse is not true in general. The example below illustrate this
fact.

\bnom\ex  {\rm  \; Let $\mathcal{X}=[1,\infty[$
 with the usual metric. Define $f:\mathcal{X}\rightarrow \mathcal{X}$ and $F:\mathcal{X}\rightarrow \mathcal{P}_{fb}(\mathcal{X})$ by, for all $x\in \mathcal{X}$,
\begin{eqnarray*}
fx=2x+1\text{\textit{, }}Fx=[1,2x+1].
\end{eqnarray*}
\begin{eqnarray*}
fx=2x+1\in Fx \text{ and } fFx=[3,4x+3] \subset Ffx=[1,4x+3].
\end{eqnarray*}
Hence, $f$ and $F$ are occasionally weakly compatible
but non weakly compatible.}
\enom

\section{General fixed point theorems}

In this section, before giving our first main result, we recall
this definition.

\bnom{\defi} Let $F:\mathcal{X} \rightarrow 2^{\mathcal{X}}$ be a set-valued map on $\mathcal{X}$. $x \in \mathcal{X}$ is a fixed point of $F$ if $x \in Fx$.
\enom
\bnom{\theo} \label{theo1} Let $f$, $g:\mathcal{X} \rightarrow
\mathcal{X}$ be maps and $F$, $G:\mathcal{X} \rightarrow
\mathcal{P}_{fb}(\mathcal{X})$ be set-valued maps such that
the pairs $\{f,F\}$ and $\{g,G\}$ are owc. Let $\varphi:(\mathbb{R^+})^5 \rightarrow
\mathbb{R}$ be a real map satisfying the following conditions:

\noindent $(\varphi _{1}):$ $\varphi$ is nonincreasing in variables $t_{4}$ and $t_{5}$,

\noindent $(\varphi _{2}):$ $\varphi(t,0,0,t,t)\geq 0$ $\forall t>0$.

\noindent If, for all $x$ and $y \in \mathcal{X}$ for which $\max\{d(fx,gy),d(fx,Fx),d(gy,Gy)\}>0 $,
$$\varphi (d(fx,gy),d(fx,Fx),d(gy,Gy),d(fx,Gy),d(gy,Fx))<0
\leqno(2.2)$$
then $f$, $g$, $F$ and $G$ have a unique common fixed point.
\enom

\noindent \demo

\noindent {\it i)} We begin to show the existence of a common fixed
point.

\noindent Since the pairs $\{f,F\}$ and $\{g,G\}$ are owc then, there exist
$u$, $v$ in $\mathcal{X}$ such that $fu\in Fu$, $gv\in Gv$,
$fFu\subseteq Ffu$ and $gGv\subseteq Ggv$.

\noindent First, we show that $gv=fu$. Suppose that is not the case, then by (2.3), we have

\begin{eqnarray*}
&&\varphi (d(fu,gv),d(fu,Fu),d(gv,Gv),d(fu,Gv),d(gv,Fu)) \\
&=&\varphi (d(fu,gv),0,0,d(fu,Gv),d(gv,Fu))<0
\end{eqnarray*}
and by $(\varphi _{1})$,
\begin{equation*}
\varphi (d(fu,gv),0,0,d(fu,gv),d(fu,gv))<0
\end{equation*}
which from $(\varphi _{2})$ gives $d(fu,gv)=0$. So $fu=gv$.

\noindent Next, we claim that
$f^{2}u=fu$. If it is not, then condition (2.3) implies that
\begin{eqnarray*}
&&\varphi (d(f^{2}u,gv),d(f^{2}u,Ffu),d(gv,Gv),d(f^{2}u,Gv),d(gv,Ffu)) \\
&=&\varphi (d(f^{2}u,fu),0,0,d(f^{2}u,Gv),d(fu,Ffu))<0.
\end{eqnarray*}
By $(\varphi _{1})$ we have
\begin{equation*}
\varphi (d(f^{2}u,fu),0,0,d(f^{2}u,fu),d(f^{2}u,fu))<0
\end{equation*}
which, from $(\varphi _{2})$,  gives $d(f^{2}u,fu)=0$. We have $f^{2}u=fu $.

\noindent Since $(f,F)$ and $(g,G)$ have the same role, we have $gv=g^2v$. Therefore, $ffu=fu=gv=ggv=gfu$, $fu=f^2u\in fFu \subset Ffu$, so
$fu\in Ffu$ and $fu=gfu \in Gfu$. Then $fu$ is a common fixed
point of $f$, $g$, $F$ and $G$.

\noindent {\it ii)} Now, we show uniqueness of the common fixed point.

\noindent Put $fu=w$ and let $w'$ be another common fixed point of the four
maps such that $w\neq w'$, then, by (2.3), we get
\begin{eqnarray*}
&&\varphi (d(fw,gw'),d(fw,Fw),d(gw',Gw'),d(fw,Gw'),d(gw',Fw)) \\
&=&\varphi (d(fw,gw'),0,0,d(fw,Gw'),d(gw',Fw))<0.
\end{eqnarray*}
By $(\varphi _{1})$, we get
\begin{equation*}
\varphi (d(fw,gw'),0,0,d(fw,gw'),d(fw,gw'))<0.
\end{equation*}
So, by ($\varphi_2$), $d(fw,gw')=0$ and thus
$d(fw,gw')=d(w,w')=0$.\qed

\bigskip

We can give two variants of Theorem \ref{theo1}:

\bnom{\theo} \label{theo2} Let $f$, $g:\mathcal{X} \rightarrow \mathcal{X}$
be maps and $F$, $G:\mathcal{X} \rightarrow \mathcal{P}_{fb}(\mathcal{X})$ be
set-valued maps such that the pairs $\{f,F\}$ and
$\{g,G\}$ are owc. Let $\varphi: (\mathbb{R^+})^6 \rightarrow \mathbb{R}$ be a real map satisfying the following conditions:

\noindent $(\varphi _{1}):$ $\varphi$ is nonincreasing in variables $t_{5}$ and $t_{6}$,

\noindent $(\varphi _{2}):$ for every $t'$, $\varphi(t',t,0,0,t,t)\geq 0$ $\forall t>0$.

\noindent If, for all $x$ and $y \in \mathcal{X}$ for which
$\max \{d(fx,gy),d(fx,Fx),d(gy,Gy)\}>0 $,
$$\varphi (H(Fx,Gy),d(fx,gy),d(fx,Fx),d(gy,Gy),d(fx,Gy),d(gy,Fx))<0
\leqno(2.3)$$
then $f$, $g$, $F$ and $G$ have a unique common fixed point.
\enom

\noindent \demo

\noindent {\it i)} We begin to show the existence of a common fixed
point in a similar proof of Theorem \ref{theo1}.

\noindent Since the pairs $\{f,F\}$ and $\{g,G\}$ are owc then, there exist
$u$, $v$ in $\mathcal{X}$ such that $fu\in Fu$, $gv\in Gv$,
$fFu\subseteq Ffu$ and $gGv\subseteq Ggv$.

\noindent First, we show that $gv=fu$. Suppose that is not the case, then condition (2.2)
implies that
\begin{eqnarray*}
&&\varphi (H(Fu,Gv),d(fu,gv),d(fu,Fu),d(gv,Gv),d(fu,Gv),d(gv,Fu)) \\
&=&\varphi (H(Fu,Gv),d(fu,gv),0,0,d(fu,Gv),d(gv,Fu))<0.
\end{eqnarray*}
By $(\varphi _{1})$ we have
\begin{equation*}
\varphi (H(Fu,Gv),d(fu,gv),0,0,d(fu,gv),d(fu,gv))<0
\end{equation*}
which from $(\varphi _{2})$ gives $d(fu,gv)=0$. So $fu=gv$.

\noindent Next, we claim that $f^{2}u=fu$. If it is not, then condition (2.2) implies that
\begin{eqnarray*}
&&\varphi (H(Ffu,Gv),d(f^{2}u,gv),d(f^{2}u,Ffu), \\
&& d(gv,Gv),d(f^{2}u,Gv),d(gv,Ffu)) \\
&=&\varphi (H(Ffu,Gv),d(f^{2}u,fu),0,0,d(f^{2}u,Gv),d(fu,Ffu))<0.
\end{eqnarray*}
By $(\varphi _{1})$ we have
\begin{equation*}
\varphi (H(Ffu,Gv),d(f^{2}u,fu),0,0,d(f^{2}u,fu),d(f^{2}u,fu))<0
\end{equation*}
which, from $(\varphi _{2})$, gives $d(f^{2}u,fu)=0$. We have $f^{2}u=fu $.

\noindent Since $(f,F)$ and $(g,G)$ have the same role, we have $gv=g^2v$.
Therefore, $ffu=fu=gv=ggv=gfu$, $fu=f^2u\in fFu \subset Ffu$, so
$fu\in Ffu$ and $fu=gfu \in Gfu$. Then $fu$ is a common fixed
point of $f$, $g$, $F$ and $G$.

\noindent {\it ii)} Now, we show uniqueness of the common fixed point.

\noindent Put $fu=w$ and let $w'$ be another common fixed point of the four maps such that $w\neq w'$, then, by (2.2), we get
\begin{eqnarray*}
&&\varphi (H(Fw,Gw'),d(fw,gw'),d(fw,Fw), \\
&&d(gw',Gw'),d(fw,Gw'),d(gw',Fw)) \\
&=&\varphi (H(Fw,Gw'),d(fw,gw'),0,0,d(fw,Gw'),d(gw',Fw))<0.
\end{eqnarray*}
By $(\varphi _{1})$, we get
\begin{equation*}
\varphi (H(Fw,Gw'),d(fw,gw'),0,0,d(fw,gw'),d(fw,gw'))<0.
\end{equation*}
So, by ($\varphi_2$), $d(fw,gw')=0$ and thus
$d(fw,gw')=d(w,w')=0$.\qed

\bnom{\theo} \label{theo3} Let $f$, $g:\mathcal{X} \rightarrow \mathcal{X}$
be maps and $F$, $G:\mathcal{X} \rightarrow \mathcal{P}_{fb}(\mathcal{X})$ be
set-valued maps such that the pairs $\{f,F\}$ and
$\{g,G\}$ are owc. Let $\varphi: (\mathbb{R^{+}})^6
\rightarrow
\mathbb{R}$ be a real map satisfying the following conditions:

\noindent $(\varphi _{1}):$ $\varphi$ is nondecreasing in variable $t_{1}$ and nonincreasing in variables $t_{5}$ and $t_{6}$,

\noindent $(\varphi _{2}):$ $\varphi(t,t,0,0,t,t)\geq 0$ $\forall\text{ \ }t>0$.

\noindent If, for all $x$ and $y \in \mathcal{X}$ for which $\max \{d(fx,gy),d(fx,Fx),d(gy,Gy)\}>0 $,
\begin{eqnarray*}
(2.4)\text{ \ }\varphi (\delta(Fx,Gy),d(fx,gy),d(fx,Fx),d(gy,Gy),d(fx,Gy),d(gy,Fx))<0
\end{eqnarray*}
then $f$, $g$, $F$ and $G$ have a unique common fixed point.
\enom

\noindent \demo

\noindent {\it i)} We begin to show existence of a common fixed point. The beginning of the proof is similar of that of previous theorems.

\noindent With the same notations, we suppose that  $gv\neq fu$.
then condition (2.4) implies that
\begin{eqnarray*}
&&\varphi (\delta(Fu,Gv),d(fu,gv),d(fu,Fu),d(gv,Gv),d(fu,Gv),d(gv,Fu)) \\
&=&\varphi (\delta(Fu,Gv),d(fu,gv),0,0,d(fu,Gv),d(gv,Fu))<0.
\end{eqnarray*}
By $(\varphi _{1})$ we have
\begin{equation*}
\varphi (d(fu,gv),d(fu,gv),0,0,d(fu,gv),d(fu,gv))<0
\end{equation*}
which from $(\varphi _{2})$ gives $d(fu,gv)=0$. So $fu=gv$. Next, we claim that
$f^{2}u=fu$. If it is not, then condition (2.4) implies that
\begin{eqnarray*}
&&\varphi (\delta(Ffu,Gv),d(f^{2}u,gv),d(f^{2}u,Ffu),d(gv,Gv),d(f^{2}u,Gv),d(gv,Ffu)) \\
&=&\varphi (\delta(Ffu,Gv),d(f^{2}u,fu),0,0,d(f^{2}u,Gv),d(fu,Ffu))<0.
\end{eqnarray*}
By $(\varphi _{1})$ we have
\begin{equation*}
\varphi (d(f^{2}u,fu),d(f^{2}u,fu),0,0,d(f^{2}u,fu),d(f^{2}u,fu))<0
\end{equation*}
which, from $(\varphi _{2})$, gives $d(f^{2}u,fu)=0$ which implies that $f^{2}u=fu $.

\noindent Since $(f,F)$ and $(g,G)$ have the same role, we have: $g^2v=gv$. Therefore, $ffu=fu=gv=ggv=gfu$, $fu=f^2u\in fFu
\subset Ffu$, so $fu\in Ffu$ and $fu=gfu \in Gfu$. Then $fu$ is a
common fixed point of $f$, $g$, $F$ and $G$.

\noindent {\it ii)} Now, we show uniqueness of the common fixed point.

\noindent Put $fu=w$ and let $w'$ be another common fixed point of the four maps such that $w\neq w'$, by (2.4), we get
\begin{eqnarray*}
&&\varphi (\delta(Fw,Gw'),d(fw,gw'),d(fw,Fw),d(gw',Gw'),d(fw,Gw'),d(gw',Fw)) \\
&=&\varphi (\delta(Fw,Gw'),d(fw,gw'),0,0,d(fw,Gw'),d(gw',Fw))<0.
\end{eqnarray*}
By $(\varphi _{1})$, we get
\begin{eqnarray*}
&&\varphi (d(fw,gw'),d(fw,gw'),0,0,d(fw,gw'),d(fw,gw')) \\
&=&\varphi (d(w,w'),d(w,w'),0,0,d(w,w'),d(w,w'))<0.
\end{eqnarray*}
So, by ($\varphi_2$), $d(w,w')=0$ and thus $w=w'$.

\bnom{\ob} {\rm Truly Theorems \ref{theo2}-\ref{theo3} are generalizations of
corresponding theorems of \cite{Ali}, \cite{Bou}, \cite{Pop1}-\cite{Pop6} and others since we extended the setting
of single-valued maps to the one of single and set-valued maps,
also we deleted the compactness in \cite{Ali}, \cite{Pop4}, we
further add that we not required the continuity, although we used
the strict contractive conditions (2.3), (2.4) which are substantially
more general than the inequalities in the cited papers, and we
weakened the concepts of compatibility, compatibility of type
$(A)$, compatibility of type $(C)$, compatibility of type $(P)$
and weak compatibility to the more general one say occasional weak compatibility. Finally we deleted some assumptions of
functions $\varphi$ which are superfluous for us but are necessary
in the papers \cite{Ali}, \cite{Bou},
\cite{Pop1}-\cite{Pop6}.}\enom

If we let $f=g$ and $F=G$ in Theorems \ref{theo1}, \ref{theo2} and \ref{theo3}, we get different corollaries.
As example, we give the following corollaries of Theorem \ref{theo3}:

\bnom{\corol} Let $f:\mathcal{X} \rightarrow \mathcal{X}$ and let
$F: \mathcal{X} \rightarrow \mathcal{P}_{fb}(\mathcal{X})$ such that the pair
$\{f,F\}$ is owc. Let $\varphi: (\mathbb{R^{+}})^6
\rightarrow \mathbb{R}$ be a real map satisfying conditions
$(\varphi _{1})$ and  $(\varphi _{2})$ of Theorem \ref{theo3} and
\begin{eqnarray*}
\varphi (\delta(Fx,Fy),d(fx,fy),d(fx,Fx),d(fy,Fy),d(fx,Fy),d(fy,Fx))<0
\end{eqnarray*}
for all $x $ and $y \in  \mathcal{X}$ for which $\max \{d(fx,fy),d(fx,Fx),d(fy,Fy)\}>0 $, then $f$ and $F$ have a
unique common fixed point.
\enom

Now, if we let $f=g$, we get the next result:

\bnom{\corol} Let $f$ be a self-map of a metric space
$(\mathcal{X},d)$ and let $F$, $G:\mathcal{X} \rightarrow
\mathcal{P}_{fb}(\mathcal{X})$ be set-valued maps. Suppose pairs $\{f,F\}$ and $\{f,G\}$
are owc and $\varphi : (\mathbb{R^{+}})^6 \rightarrow \mathbb{R}$
satisfies conditions $(\varphi _{1})$ and $(\varphi _{2})$ of Theorem \ref{theo3} and
\begin{eqnarray*}
\varphi (\delta(Fx,Gy),d(fx,fy),d(fx,Fx),d(fy,Gy),d(fx,Gy),d(fy,Fx))<0
\end{eqnarray*}
for all $x$ and $y \in \mathcal{X}$ for which $\max \{d(fx,fy),d(fx,Fx),d(fy,Gy)\}>0$, then $f$, $F$ and $G$ have
a unique common fixed point.
\enom

With different choices of the real map $\varphi$, we obtain the
following corollaries:

\bnom{\corol}\label{corol1} If in the hypotheses of Theorem
\ref{theo3}, we have instead of (2.4) one of the following
inequalities, for all $x$ and $y \in \mathcal{X}$ whenever the right hand side of
each inequality is not zero, then the four maps have a unique common fixed point.
\begin{eqnarray*}
(a)\text{ \ }\delta(Fx,Gy)< k \max \{d(fx,gy),d(fx,Fx),d(gy,Gy),d(fx,Gy),d(gy,Fx)\}
\end{eqnarray*}
where $0<k\leq 1$,
\begin{eqnarray*}
(b)\text{ \ }\delta^{2}(Fx,Gy) &<& c_{1}\max\{d^{2}(fx,gy),d^{2}(fx,Fx),d^{2}(gy,Gy)\} \\
&&+c_{2}\max\{d(fx,Fx)d(fx,Gy),d(gy,Gy)d(gy,Fx)\} \\
&&+c_{3}d(fx,Gy)d(gy,Fx)
\end{eqnarray*}
where $c_{1}>0$, $c_{2}$, $c_{3}\geq 0$ and $c_{1}+c_{3}\leq 1$,
\begin{eqnarray*}
(c)\text{ \ }\delta(Fx,Gy)&<&[ \alpha \delta^{p-1}(Fx,Gy)d(fx,gy) \\
&&+\beta \delta^{p-2}(Fx,Gy)d(fx,Fx)d(gy,Gy)\\
&&+\gamma d^{p-1}(fx,Gy)d(gy,Fx) \\
&&+\delta d(fx,Gy)d^{p-1}(gy,Fx)]^{\frac{1}{p}}
\end{eqnarray*}
where $\alpha>0$, $\beta$, $\gamma$, $\delta\geq 0$, $\alpha+\gamma+\delta\leq 1$ and $p\geq 2$,
\begin{eqnarray*}
(d)\text{ \ }\delta^{2}(Fx,Gy) <\frac{1}{\alpha}\left[\beta d^{2}(fx,gy)
+\frac{\gamma d(fx,Gy)d(gy,Fx)}{1+\delta d^{2}(fx,Fx)+\epsilon
d^{2}(gy,Gy)}\right]
\end{eqnarray*}
where $\alpha>0$, $\beta,\gamma,\delta,\epsilon\geq 0$ and $\beta+\gamma\leq \alpha$,
\begin{eqnarray*}
(e)\text{ \ }\delta(Fx,Gy) &<&\left[\alpha d^{p}(fx,gy)+\beta d^{p}(fx,Fx)+\gamma d^{p}(gy,Gy)
\right] ^{\frac{1}{p}} \\
&&+\delta \left[d(fx,Gy)d(gy,Fx)\right]^{\frac{1}{2}}
\end{eqnarray*}
where $0<\alpha \leq (1-\delta)^{p}$, $\beta,\;\gamma,\;
\delta\geq 0$ and $ p\in {\mathbb{N}}^{\ast }=\{1,2,\ldots\}.$
\enom

\noindent \demo

\noindent For proof of $(a)$, $(b)$, $(c)$, $(d)$ and $(e)$, we use Theorem \ref{theo3} with the following
functions $\varphi$ which satisfy, for every case,
hypothesis $(\varphi _{1})$ and $(\varphi _{2})$

\noindent for $(a)$:
\begin{eqnarray*}
&&\varphi (\delta(Fx,Gy),d(fx,gy),d(fx,Fx),d(gy,Gy),d(fx,Gy),d(gy,Fx)) \\
&=&\delta(Fx,Gy)-k\max \{d(fx,gy),d(fx,Fx),d(gy,Gy),d(fx,Gy),d(gy,Fx)\}.
\end{eqnarray*}
This function $\varphi$ is used by many authors with single
maps, for example: \cite {JuR2} in Theorem 1, Example 3.4 in \cite {PRV}.

\noindent For $(b)$:
\begin{eqnarray*}
&&\varphi (\delta(Fx,Gy),d(fx,gy),d(fx,Fx),d(gy,Gy),d(fx,Gy),d(gy,Fx)) \\
&=&\delta^{2}(Fx,Gy)-c_{1}\max \{d^{2}(fx,gy),d^{2}(fx,Fx),d^{2}(gy,Gy)\} \\
&&-c_{2}\max \{d(fx,Fx)d(fx,Gy),d(gy,Gy)d(gy,Fx)\} \\
&&-c_{3}d(fx,Gy)d(gy,Fx).
\end{eqnarray*}
This function $\varphi$ is Example 2 of \cite{Pop4}.

\noindent For $(c)$:
\begin{eqnarray*}
&&\varphi (\delta(Fx,Gy),d(fx,gy),d(fx,Fx),d(gy,Gy),d(fx,Gy),d(gy,Fx)) \\
&=&\delta(Fx,Gy)-\left[\alpha \delta^{p-1}(Fx,Gy)d(fx,gy)\right. \\
&&+\beta \delta^{p-2}(Fx,Gy)d(fx,Fx)d(gy,Gy) \\
&&\left. +\gamma d^{p-1}(fx,Gy)d(gy,Fx)+\delta d(fx,Gy)d^{p-1}(gy,Fx)\right] ^{\frac{1}{p}}.
\end{eqnarray*}
For $p=3$, we have Example 3.4 of \cite{Bou} and Example 3 of \cite {Pop5}. If we take $p=2$, $\varphi$ is Example 1 of
\cite{Pop2}.

\noindent For $(d)$:
\begin{eqnarray*}
&&\varphi (\delta(Fx,Gy),d(fx,gy),d(fx,Fx),d(gy,Gy),d(fx,Gy),d(gy,Fx)) \\
&=&\delta^{2}(Fx,Gy)-\frac{1}{\alpha}\left[\beta d^{2}(fx,gy)+\frac{\gamma d(fx,Gy)d(gy,Fx)}{1+\delta d^{2}(fx,Fx)+\epsilon d^{2}(gy,Gy)}\right].
\end{eqnarray*}
This function $\varphi$ is that one of Example 6 of \cite{Pop1}.

\noindent And for $(e)$:
\begin{eqnarray*}
&&\varphi (\delta(Fx,Gy),d(fx,gy),d(fx,Fx),d(gy,Gy),d(fx,Gy),d(gy,Fx)) \\
&=&\delta(Fx,Gy)-\left[\alpha d^{p}(fx,gy)+\beta d^{p}(fx,Fx)+\gamma d^{p}(gy,Gy)\right] ^{\frac{1}{p}} \\
&&-\delta \left[d(fx,Gy)d(gy,Fx)\right] ^{\frac{1}{2}}
\end{eqnarray*}
\qed

\bnom{\corol} Let $f$, $g$ be two self-maps of a
metric space $(\mathcal{X},d)$ and let $F$ and $G:\mathcal{X}
\rightarrow \mathcal{P}_{fb}(\mathcal{X})$ be set-valued maps such that the
pairs $\{f,F\}$ and $\{g,G\}$ are owc.
Suppose that, for all $x$, $y \in \mathcal{X}$, we have the
inequality
\begin{eqnarray*}
(f)\text{ \ }\delta^{p}(Fx,Gy)<\alpha d^{p}(fx,gy)+\beta d^{p}(fx,Fx)+\gamma d^{p}(gy,Gy)
\end{eqnarray*}
such that $0<\alpha \leq 1$, $\beta$ and $ \gamma\geq 0$ and $p\in \mathbb{N}^{\ast }=\{1,2,\ldots\}$ whenever the right hand side of
the above inequality is positive. Then $f$, $g$, $F$ and $G$ have
a unique common fixed point.
\enom

\noindent \demo

\noindent We give this corollary because it is an interesting particular case
of the previous corollary. We obtain the result by using (e) in
Corollary \ref{corol1} with $\delta=0$.\qed

\section{Two other type common fixed point theorems}

We begin by a Gregu\v{s} type common fixed point theorem. As we already said, in 1998, Pathak et al. \cite{PCKM} gave an extension of
compatibility of type $(A)$ by introducing the concept of
compatibility of type $(C)$ and they proved a common fixed point
theorem of Gregu\v{s} type for four compatible maps of type $(C)$
in a Banach space. Further, Djoudi and Nisse \cite{DjN} extended
the result of \cite{PCKM} by weakening compatibility of type $(C)$
to the weak one without continuity.

Our objective here is to establish a common fixed point theorem
for four occasionally weakly compatible single and set-valued
maps of Gregu\v{s} type in a metric space which improves
the results of \cite{DjN}, \cite{PCKM} and others.

\bnom{\theo} \label{theo4} Let $f$ and $g : \mathcal{X}
\rightarrow \mathcal{X}$ be maps, $F$ and $G : \mathcal{X}
\rightarrow \mathcal{P}_{fb}(\mathcal{X})$ be set-valued maps such that the
pairs $\{f,F\}$ and $\{g,G\}$ are owc.
Let $\Psi:\mathbb{R}^{+}\rightarrow \mathbb{R}^{+}$ be a nondecreasing map such that, for every $t>0$, $\Psi (t)<t$ and
satisfying the following condition:
$$
\delta^{p}(Fx,Gy) \leq
\Psi [ad^{p}(fx,gy)+(1-a)d^{\frac{p}{2}}(gy,Fx)d^{\frac{p}{2}}(fx,Gy)]
\leqno(3.1)$$
for all $x$ and $y \in \mathcal{X}$, where $0<a\leq 1$
and $p\geq 1$. Then $f$, $g$, $F$ and $G$ have a unique common fixed point.
\enom

\noindent \demo

\noindent Since $f$, $F$ and  $g$, $G$  are owc, as in proof of Theorem \ref{theo2}, there exist $u$, $v$ in $\mathcal{X}$ such that $fu\in Fu$, $gv\in Gv$, $fFu\subseteq Ffu$, $gGv\subseteq Ggv$.

\noindent {\it i)} As in proof of Theorem \ref{theo2}, we begin to show existence of a common fixed
point. We have,
\begin{eqnarray*}
\delta^p(Fu,Gv)\leq \Psi (ad^p(fu,gv)+(1-a)d^{\frac{p}{2}}(gv,Fu) d^{\frac{p}{2}}(fu, Gv))
\end{eqnarray*}
and by the properties of $\delta$ and $\Psi$, we get
\begin{eqnarray*}
d^p(fu,gv)\leq \delta^p(Fu,Gv)\leq \Psi (d^p(fu,gv)).
\end{eqnarray*}
So, if $d(fu,gv)>0$,  $\Psi (t)<t$ for $t>0$, we obtain
\begin{eqnarray*}
d^p(fu,gv)\leq \delta^p(Fu,Gv)\leq \Psi (d^p(fu,gv))<d^p(fu,gv)
\end{eqnarray*}
which is a contradiction, thus we have $d(fu,gv)=0$, hence $fu=gv$.

\noindent Again, if $d(f^2u,fu)>0$, then by (3.1), we have
\begin{eqnarray*}
\delta^p(Ffu,Gv) \leq \Psi [a d^p(f^2u, gv)+(1-a) d^{\frac{p}{2}}(gv,Ffu) d^{\frac{p}{2}}(f^2u,Gv)]
\end{eqnarray*}
and hence
\begin{eqnarray*}
d^p(f^2u,fu)\leq \delta^p(Ffu,Gv) \leq \Psi [d^p(f^2u,fu)]
\end{eqnarray*}
Since $d(f^2u,fu)>0$, we obtain
\begin{eqnarray*}
d^p(f^2u,fu)\leq \delta^p(Ffu,Gv) \leq \Psi [d^p(f^2u,fu)]<d^p(f^2u,fu)
\end{eqnarray*}
what it is impossible. Then we have $d(f^2u,fu)=0$; i.e., $f^2u=fu$. Similarly, we can prove that $g^{2}v=gv$, let $fu=w$ then, $fw=w=gw$, $w\in Fw$ and $w\in Gw$, this completes the proof of the existence.

\noindent {\it ii)} For the uniqueness, let $w'$ be a second common fixed point of $f$, $g$, $F$ and $G$ with $w'\neq w$. Then, $d(w,w')=d(fw,gw')\leq \delta(Fw,Gw')$ and, by assumption (3.1), we obtain
\begin{eqnarray*}
\delta^{p}(Fw,Gw') \leq \Psi [ad^{p}(fw,gw')+(1-a)d^{\frac{p}{2}}(fw,Gw') d^{\frac{p}{2}}(gw',Fw)]
\end{eqnarray*}
and thus
\begin{eqnarray*}
d^{p}(w,w')=d^{p}(fw,gw')\leq \delta^{p}(Fw,Gw') \leq \Psi [d^{p}(w,w')]<d^{p}(w,w')
\end{eqnarray*}
Since $d(w,w')>0$, we have a contradiction. So, $w=w'$.\qed

\bnom{\theo} \label{theo5} Let $f$ and $g:\mathcal{X} \rightarrow \mathcal{X}$ be maps, $F$ and $G:\mathcal{X} \rightarrow \mathcal{P}_{fb}(\mathcal{X})$ be set-valued maps
such that the pairs $\{f,F\}$ and $\{g,G\}$ are owc. Let $\Psi:\mathbb{R}^{+}\rightarrow
\mathbb{R}^{+}$ be a nondecreasing map such that, for every
$t>0$, $\Psi (t)<t$ and satisfying the following condition:
\begin{eqnarray*}
\delta^{p}(Fx,Gy) &\leq& \Psi [ad^{p}(fx,gy)+(1-a)\max\{\alpha d^{p}(fx,Fx),\beta d^{p}(gy,Gy), \\
&&d^{\frac{p}{2}}(fx,Fx)d^{\frac{p}{2}}(gy,Fx),d^{\frac{p}{2}}(gy,Fx)d^{\frac{p}{2}}(fx,Gy), \\
&&\frac{1}{2}(d^{p}(fx,Fx)+d^{p}(gy,Gy))\}]
\end{eqnarray*}
for all $x$ and $y \in \mathcal{X}$, where $0<a\leq 1$, $0<\alpha ,\beta \leq 1$ and $p\geq 1$. Then $f$, $g$, $F$ and $G$ have a unique common fixed point.
\enom

\noindent \demo

\noindent Since $f$, $F$ and  $g$, $G$  are owc, as in proof of Theorem \ref{theo1}, there exist $u$, $v$ in $\mathcal{X}$ such that $fu\in Fu$, $gv\in Gv$, $fFu\subseteq Ffu$, $gGv\subseteq Ggv$. Since $\Psi $ is a nondecreasing function and since for
any real numbers $c$ and $d$, $\frac{c+d}{2}\leq \max \{c,d\}$ we have, for all $x$, $y \in \mathcal{X}$,
\begin{eqnarray*}
\delta^{p}(Fx,Gy) &\leq& \Psi [ad^{p}(fx,gy)+(1-a)\max \{d^{p}(fx,Fx),d^{p}(gy,Gy) \\
&&d^{\frac{p}{2}}(fx,Fx)d^{\frac{p}{2}}(gy,Fx),d^{\frac{p}{2}}(gy,Fx)d^{\frac{p}{2}}(fx,Gy) \}]
\end{eqnarray*}
and, for $u$ and $v$,
\begin{eqnarray*}
\delta^{p}(Fu,Gv) \leq \Psi [ad^{p}(fu,gv)+(1-a)d^{\frac{p}{2}}(gv,Fu)d^{\frac{p}{2}}(fu,Gv)].
\end{eqnarray*}
The continuation of the proof is identical of that of Theorem \ref{theo2}. \qed
\bigskip

 If in (3.1), we replace $\delta$ with $H$ and $\Psi (t)< t$ with $\Psi (t)\leq t$, we can prove the
\bnom{\theo} Let $f$ and $g : \mathcal{X} \rightarrow \mathcal{X}$
be maps, $F$ and $G : \mathcal{X} \rightarrow
\mathcal{P}_{fb}(\mathcal{X})$ be set-valued maps such that the
pairs $\{f,F\}$ and $\{g,G\}$ are owc. Let  $u$ and $v$ in
$\mathcal{X}$ such that $fu\in Fu$, $gv\in Gv$, $fFu\subseteq
Ffu$, $gGv\subseteq Ggv$. Let $\Psi:\mathbb{R}^{+}\rightarrow
\mathbb{R}^{+}$ be a nondecreasing map such that, for every $t>0$,
$\Psi (t)\leq t$ and satisfying the following condition:
$$ H^{p}(Fx,Gy) \leq \Psi
[ad^{p}(fx,gy)+(1-a)d^{\frac{p}{2}}(gy,Fx)d^{\frac{p}{2}}(fx,Gy)]
\leqno(3.3) $$
 for all $x$ and $y \in \mathcal{X}$, where $0<a\leq
1$ and $p\geq 1$.\\ If $fu=gv$ is a common fixed point of $f$ and
$g$, then $fu$ is a common fixed point of $f$, $g$, $F$ and $G$
and $Fu=Gv$. \enom
\demo

\noindent Since $gv\in Gv$, $fu\in Fu$ and $f^2u \in fFu \subseteq Ffu$, we have \\
$d(gv, Fu)\leq H(Fu, Gv)$, $d(fu, Gv)\leq H(Fu, Gv)$,\\
$d(gv, Ffu)\leq H(Ffu, Gv)$ and $ d(f^2u, Gv) \leq H(Ffu, Gv)$.\\
From the nondecrease of $\Psi$, we obtain
\begin{eqnarray*}
H^p(Ffu,Gv) &\leq &\Psi[a d^p(f^2u, gv)+(1-a)d^{\frac{p}{2}}(gv,Ffu)d^{\frac{p}{2}}(f^2u,Gv)]\\
&\leq& \Psi[a d^p(f^2u, gv)+(1-a)H^p(Ffu,Gv)]
\end{eqnarray*}

$$H^p(Fu,Gv) \leq \Psi[a d^p(fu, gv)+(1-a)H^p(Fu,Gv)]$$
and
$$ H^p(Fu,Ggv) \leq \Psi[a d^p(fu, g^2v)+(1-a)H^p(Fu,Ggv)].$$
Now, if $Fu\neq Gv$, since, for every $t>0$,
$\Psi (t)\leq t$, $$H^p(Fu,Gv) \leq a d^p(fu, gv)+(1-a)H^p(Fu,Gv).$$
Consequently, $H(Fu,Gv) \leq d(fu,gv)$ and $fu \neq gv$. We have shown that if $fu=gv$, then $Fu=Gv$.
By similar proofs, if $f^2u=gv$,  then $Gv=Ffu$ and if $fu=g^2v$, then $Fu= Ggv$. The proof is finished.\qed

\bnom{\ob} {\rm Obviously, Theorems \ref{theo4} and \ref{theo5}
extend the results of \cite{DjN}, \cite{PCKM} and others to the
class of four single and set-valued maps. In particular, Theorem
\ref{theo5} improves the cited results since we not required the
closeness of the sets $F(\mathcal{X})$ and $G(\mathcal{X})$, also
we deleted the inclusions $F(\mathcal{X})\subset f(\mathcal{X})$
and $G(\mathcal{X})\subset g(\mathcal{X})$ in \cite{DjN}, we
weakened the weak compatibility in \cite{DjN} and the
compatibility of type (C) in \cite{PCKM} to the wider one cited
occasional weak compatibility and we deleted the continuity
which is indispensable in \cite{PCKM} and the upper semicontinuity
imposed on $\Psi $ in \cite{DjN}.} \enom

If we put $f=g$ in Theorem \ref{theo4}, then we get the corollary:

\bnom{\corol} Let $f:\mathcal{X} \rightarrow \mathcal{X}$ be a map and let $F$ and $G:\mathcal{X} \rightarrow
\mathcal{P}_{fb}(\mathcal{X})$ be set-valued maps. Let
$\Psi:\mathbb{R}^{+}\rightarrow \mathbb{R}^{+}$ be a nondecreasing
map such that, for every $t>0$, $\Psi (t)<t$. Suppose pairs
$\{f,F\}$ and $\{f,G\}$ are owc and
satisfy the inequality
\begin{eqnarray*}
\delta^{p}(Fx,Gy) \leq \Psi [ad^{p}(fx,fy)+(1-a)d^{\frac{p}{2}}(fy,Fx)d^{\frac{p}{2}}(fx,Gy)]
\end{eqnarray*}
for all $x$, $y \in \mathcal{X}$, where
$0<a\leq 1$ and $p\geq 1$, then $f$, $F$ and $G$ have a unique common fixed point.
\enom

If we put $f=g$ and $F=G$ in Theorem \ref{theo4}, then we obtain
the following result:

\bnom{\corol} Let $f:\mathcal{X} \rightarrow \mathcal{X}$ be a map and let $F:\mathcal{X} \rightarrow
\mathcal{P}_{fb}(\mathcal{X})$ be set-valued mapping such that $f$ and $F$ are
owc . Let $\Psi:\mathbb{R}^{+}\rightarrow \mathbb{R}^{+}$ be a
nondecreasing map such that, for every $t>0$, $\Psi (t)<t$. If
\begin{eqnarray*}
\delta^{p}(Fx,Fy) &\leq &\Psi [ad^{p}(fx,fy)+(1-a)\max \{\alpha
d^{p}(fx,Fx),\beta d^{p}(fy,Fy), \\
&&d^{\frac{p}{2}}(fx,Fx)d^{\frac{p}{2}}(fy,Fx),d^{\frac{p}{2}}(fy,Fx)d^{\frac{p}{2}}(fx,Fy), \\
&& \frac{1}{2}(d^{p}(fx,Fx)+d^{p}(fy,Fy))\}]
\end{eqnarray*}
for all $x$, $y \in \mathcal{X}$, where $0<a\leq 1$, $\{\alpha ,\beta \} \subset ]0,1]$ and $p\geq 1$, then $f$ and $F$ have a
unique common fixed point.
\enom

Now, we end our work by establishing a near-contractive common fixed point theorem which improves
those given by Elamrani and Mehdaoui \cite{ElM}, Mbarki \cite{Mba}
and others since our version does not impose continuity and we use
occasional weak compatibility which is more general than
compatibility and weak compatibility; also we delete, on $\Phi$,
some strong conditions which are necessary in papers \cite{ElM}
and \cite{Mba} on a metric space instead of a complete metric
space.

\bnom\theo  \label{theo6} Let $f$ and $g:\mathcal{X} \rightarrow \mathcal{X}$
be maps, $F$ and $G:\mathcal{X} \rightarrow
\mathcal{P}_{fb}(\mathcal{X})$ be set-valued maps and $\Phi$ be a
nondecreasing function of $[0,\infty[$ into itself such that
$\Phi(t)=0$ iff $t=0$ and satisfying inequality
\begin{eqnarray*}
(3.7)\text{ \ }\; \; \; \Phi (\delta(Fx,Gy)) &\leq& \alpha (d(fx,gy))\Phi (d(fx,gy)) \\
&&+\beta(d(fx,gy))[\Phi (d(fx,Gy))+\Phi (d(gy,Gy))] \\
&&+\gamma(d(fx,gy))[\Phi (d(fx,Fx))+\Phi (d(gy,Fx))]
\end{eqnarray*}
for all $x$, $y \in \mathcal{X}$ and $\alpha$, $\beta$, $\gamma:[0,\infty[ \rightarrow [0,1[$ satisfying condition
\begin{eqnarray*}
(4)\text{ \ }\alpha(t)+\beta(t)+\gamma(t)<1\text{ \ }\forall t>0.
\end{eqnarray*}
If the pairs $\{f,F\}$ and $\{g,G\}$ are owc, then $f$, $g$, $F$ and $G$ have a unique common fixed point in $\mathcal{X}$.
\enom

\noindent \demo

\noindent Since $f$, $F$ and  $g$, $G$  are owc, as in proof of Theorem \ref{theo1}, there exist $u$, $v$ in $\mathcal{X}$ such that $fu\in Fu$, $gv\in Gv$, $fFu\subseteq Ffu$, $gGv\subseteq Ggv$.

\noindent {\it i)} First we prove that $fu=gv$. By (3.7), we have
\begin{eqnarray*}
\Phi (\delta(Fu,Gv)) &\leq &\alpha(d(fu,gv))\Phi (d(fu,gv)) \\
&&+\beta(d(fu,gv))[\Phi (d(fu,Gv))+\Phi (d(gv,Gv))]  \\
&&+\gamma(d(fu,gv))[\Phi (d(fu,Fu))+\Phi (d(gv,Fu))] \\
&=&\alpha(d(fu,gv))\Phi (d(fu,gv))+\beta(d(fu,gv))\Phi (d(fu,Gv)) \\
&&+\gamma(d(fu,gv))\Phi (d(gv,Fu)).
\end{eqnarray*}
If $d(fu,gv)>0 $, since $\Phi$ is nondecreasing and $\Phi(t)=0$ iff $t=0$, from inequalities (3.7) and (4) we get
\begin{eqnarray*}
\Phi(d(fu,gv))&\leq& \Phi (\delta(Fu,Gv)) \\
&\leq&\alpha(d(fu,gv))\Phi (d(fu,gv))+\beta(d(fu,gv))\Phi (d(fu,Gv)) \\
&&+\gamma(d(fu,gv))\Phi (d(gv,Fu)) \\
&\leq& [\alpha(d(fu,gv))+\beta(d(fu,gv))+\gamma(d(fu,gv))] \Phi (d(fu,gv)) \\
&<&\Phi(d(fu,gv))
\end{eqnarray*}
which is a contradiction. Hence $d(fu,gv)=0$ and thus $fu=gv$.

\noindent Now we claim that $f^{2}u=fu$. Suppose not, since $\Phi $ is nondecreasing and $\Phi(t)=0$
iff $t=0$, the use of (3.7) and (4) gives
\begin{eqnarray*}
\Phi(d(f^2u,fu))&\leq& \Phi (\delta(Ffu,Gv)) \\
&\leq& \alpha(d(f^{2}u,gv))\Phi (d(f^{2}u,gv)) \\
&&+\beta(d(f^{2}u,gv))[\Phi (d(f^{2}u,Gv))+\Phi (d(gv,Gv))] \\
&&+\gamma(d(f^{2}u,gv))[\Phi (d(f^{2}u,Ffu))+\Phi (d(gv,Ffu))] \\
&=&\alpha(d(f^{2}u,fu))\Phi (d(f^{2}u,fu))+\beta(d(f^{2}u,fu))\Phi (d(f^{2}u,Gv)) \\
&&+\gamma(d(f^{2}u,fu))\Phi (d(fu,Ffu)) \\
&\leq &[\alpha(d(f^{2}u,fu))+\beta(d(f^{2}u,fu)) \\
&&+\gamma(d(f^{2}u,fu))] \Phi (d(f^{2}u,fu)) \\
&<& \Phi(d(f^2u,fu))
\end{eqnarray*}
this contradiction implies that $\Phi(d(f^2u,fu))=0$ and hence
$f^{2}u=fu$. Similarly, we can prove that $g^{2}v=gv$. So, if
$w=fu=gv$ therefore $fw=w=gw$, $w\in Fw$ and $w\in Gw$. Existence
of a common fixed point is proved.

\noindent {\it ii)} Assume that there exists a second common fixed point
$w'$ of $f$, $g$, $F$ and $G$ such that $w'\neq w$. We have
$d(w,w')=d(fw,gw')\leq \delta(Fw,Gw')$. Since $d(w,w')>0$, by
inequality (3.7) and properties of functions $\Phi$, $\alpha$ and $\beta$, we obtain
\begin{eqnarray*}
\Phi(d(w,w'))&\leq& \Phi (\delta(Fw,Gw')) \\
&\leq& \alpha(d(fw,gw'))\Phi (d(fw,gw')) \\
&&+\beta(d(fw,gw'))[\Phi (d(fw,Gw'))+\Phi (d(gw',Gw'))] \\
&&+\gamma(d(fw,gw'))[\Phi (d(fw,Fw))+\Phi (d(gw',Fw))] \\
&=&\alpha(d(w,w'))\Phi (d(w,w'))+\beta(d(w,w'))\Phi (d(w,Gw')) \\
&&+\gamma(d(w,w'))\Phi (d(w',Fw)) \\
&\leq& [\alpha(d(w,w'))+\beta(d(w,w'))+\gamma(d(w,w'))] \Phi (d(w,w')) \\
&<& \Phi(d(w,w'))
\end{eqnarray*}
this contradiction implies that $\Phi(d(w,w'))=0$, hence $w'=w$.\qed

\bnom{\ob } {\rm  The above theorem remains valid if we replace inequality (3.7) by the following one:
\begin{eqnarray*}
\Phi (\delta(Fx,Gy)) &\leq& \alpha (d(fx,gy))\Phi (d(fx,gy)) \\
&&+\beta (d(fx,gy))\max\{\Phi (d(fx,Gy)),\Phi (d(gy,Gy))\} \\
&&+\gamma (d(fx,gy))[\Phi (d(fx,Fx))+\Phi (d(gy,Fx))].
\end{eqnarray*}}
\enom

\end{document}